\documentclass[reqno,a4paper,12pt]{amsart}
\usepackage{color}
\usepackage[latin1]{inputenc}
\usepackage[english]{babel}
\usepackage{amsmath,amssymb,amsfonts,amsthm,enumerate}
\usepackage{mathrsfs}
\usepackage[T1]{fontenc}
\usepackage[active]{srcltx}
\usepackage{epsfig}
\usepackage{graphicx}

\numberwithin{equation}{section} 
\newtheorem{theorem}{Theorem}[section]

\newtheorem{example}[theorem]{Example}

\newcommand{\R}{\mathbb{R}}

 \newcommand{\tauV}{{\kern-3pt\tau}}

 \newcommand{\oVVVk}{\overline{\mbox{\boldmath$V$}}\kern-3pt}
 \newcommand{\tVVVk}{\tilde{\mbox{\boldmath$V$}}\kern-3pt}

 \newcommand\mpar[1]{\marginpar{\tiny \color{red} #1}}
 \newcommand\bblue[1]{{\color{blue} #1}}

\begin{document}

\title[Unconstrained Free boundary problems]{An overview of unconstrained\\ free boundary problems}

\author[Alessio Figalli]{Alessio Figalli}
\address{Mathematics Department, The University of Texas at Austin,  Austin, Texas, 78712-1202, USA}
\email{figalli@math.utexas.edu}

\author[Henrik Shahgholian ]{Henrik Shahgholian}
\address{Department of Mathematics, KTH Royal Institute of Technology, 100~44  Stockholm, Sweden}
\email{henriksh@kth.se}

\thanks{A. Figalli was partially supported by the NSF grants DMS-1262411 and DMS-1361122.
H. Shahgholian was partially supported by the Swedish Research Council.}

\begin{abstract}    
In this paper we present a survey  concerning  unconstrained free boundary problems of type
 $$ \left\{
\begin{array}{ll}
F_1(D^2u,\nabla u,u,x)=0 & \text{in }B_1 \cap \Omega ,\\
F_2 (D^2 u,\nabla u,u,x)=0 & \text{in }B_1\setminus\Omega ,\\
u \in \mathbb{S}(B_1),
\end{array}
\right.
$$
where $B_1$ is the unit ball,   $\Omega$ is an unknown open set, $F_1, F_2$ are elliptic operators (admitting regular solutions),
and $\mathbb{S}$ is a functions space to be specified in each case. 
Our main objective is to discuss a unifying approach to the optimal regularity of solutions to the above matching problems,
and list several open problems in this direction. 
\end{abstract}

\maketitle


\section{Introduction}


\subsection{Background}

Free boundary problems arise naturally in a number of physical phenomena. The common theme in these problems is to find  an unknown   pair $(u,\Omega)$, where the function $u$ solves some equation outside $\partial \Omega$ (the free boundary), 
and some global conditions determine the behavior of $u$ across $\partial \Omega$.

Most of these problems fit into two categories: the class of obstacle-type problems, and the Bernoulli-type ones.
It is worth noticing that there are several other free boundary problems, such as free boundaries in porous medium equations, curvature related problems, and so on. However in this note we focus mainly on obstacle-type problems, and just mention in passing the Bernoulli-type ones.


\subsection{Obstacle problems} 
The prototype model for the first class of problems is given by the classical obstacle problem
\begin{equation}\label{eq:obst}
\min_{\mathcal K} \int_{B_1} |\nabla v|^2 \qquad  {\mathcal K}:=\{v \in W^{1,2}(B_1) :\  v \geq \psi,\,v|_{\partial B_1}=g\},
\end{equation}
where $B_1$ is the unit ball,  $\psi:B_1\to {\mathbb{R}}$ and $g:\partial B_1\to {\mathbb{R}}$ are smooth given functions, and $\psi|_{\partial B_1}<g$. This class of problems is 
very well described in the book \cite{PSU}, and we refer to it for more details and references.

It is well known that if $u$ is the unique minimizer to the above problem, then it solves 
$$
\Delta u=0 \qquad \text{in } \{u>\psi\} \cap B_1.
$$
It is possible to prove that solutions to this problem belong to $W^{2,p} (B_1)$ for any $p>1$,
and because $W^{2,p}$ functions are twice  differentiable a.e. when $p> n$, one deduces that 
$\Delta u=\Delta \psi$ a.e. on $\{u=\psi\}$.
Hence the obstacle problem can be rewritten as
$$
\Delta u=\Delta \psi\,\chi_{\{u=\psi\}}  \qquad \text{in } B_1.
$$
with the additional (global) information $u \in W^{2,p}(B_1)$ for all $p\in (1,\infty)$ \cite[Chapter 1, Section 3]{PSU}.
Notice that because $\Delta u$ jumps from the value zero inside $\{u>\psi\}$ to the value $\Delta\psi$ on $\{u=\psi\}$,
the second derivatives of $u$ cannot be continuous in general.
However, it is possible to prove that $u \in C^{1,1}$, and starting from this fact one can
show the regularity of the free boundary $\partial \{u>\psi\}$  \cite{Caff,Caobstacle} (see also \cite[Chapters 2-8]{PSU}).\\

\subsection{Bernoulli problems} 
Bernoulli-type problems are free boundary problems with a transition conditions across the free boundary.
 The archetype for this class of problems is 
\begin{equation}\label{eq:obst}
\min_{\mathcal K} \int_{B_1} |\nabla v|^2+\chi_{\{v>0\}} \qquad  {\mathcal K}:=\{v \in W^{1,2}(B_1) :\  v|_{\partial B_1}=g\},
\end{equation}
where $g$ may take both negative and positive values.
Solutions to this problem satisfy
$$
\Delta u=0\qquad \text{when $u \neq 0$},
$$
with a transition condition  on the free boundary $\{u=0\}$ that can be derived formally in the case of a smooth free boundary and reads
$$
(u_\nu^+)^2-(u_\nu^-)^2=1.
$$
Here $u_\nu^\pm$  denote  the normal derivatives in the inward direction to $\{\pm u>0\}$,
so that $u_\nu^\pm$  are both nonnegative. 

Because of the jump in the gradient along the free boundary, the optimal possible regularity is Lipschitz.\footnote{It is worth noticing  that classical PDE theory usually deals with regularity of type $W^{k,p}$ or $C^{k,\alpha}$ with $1<p<\infty$ and $\alpha \in (0,1)$.
On the other hand, free boundary problems  give rise to  integer-order regularity ($C^{0,1}$ in this case, or $C^{1,1}$ in the obstacle problem), and as such these regularities are much harder to obtain, since classical techniques usually fail in such scenarios.}
This regularity is indeed true \cite[Chapter 6]{CSbook}, and again regularity of the free boundary can be shown
\cite[Chapters 3-5]{CSbook}. We refer  to the book \cite{CSbook} for more details and references.\\

\subsection{Obstacle-type problems}
The aim of this note is to focus on the first class of problems, although an analysis (rather a discussion, since the second class of problems seems much harder to handle)  similar to the one performed here could be done also for the second class. 
The basic idea is that solutions to the classical obstacle problem solve the system
$$
\left\{
\begin{array}{ll}
\Delta u=0 & \text{in }B_1 \cap \Omega ,\\
u=\psi & \text{in }B_1\setminus\Omega ,\\
u \in W^{2,p}(B_1)&  \forall\,p\in (1,\infty),
\end{array}
\right.
$$
where $\Omega=\{u>\psi\}$. However, in the system above we can neglect the information that $\{u>\psi\}$ inside $\Omega$
and just ask ourselves what is the regularity of $u$ if $\Omega$ is some arbitrary open set (in particular,
$u$ could be less than $\psi$ in some regions of $\Omega$).
This problem is related to free boundaries in potential theory and as shown in the discussion of Problem A in \cite{PSU}, even if we are losing some information,
under the assumption that $\psi\in C^{1,1}$
one can still prove that solutions to this problem are $C^{1,1}$,
and an analysis of the free boundary can still be performed.

One may now try to go further and replace the Laplacian by a fully nonlinear operator. Actually, also the condition $u=\psi$
can be seen as a degenerate PDE, as we shall see below.\\

Our goal here is to present a general class of free boundary problems as a matching problem.
This  matching   is across an unknown (free) boundary,
in the sense that one is looking for a function which satisfies two different conditions in disjoint domains
but with some global information about how the values of the function in the two domains should match.
Such matching problems can be represented in a very general form as
\begin{equation}\label{matching}
\left\{
\begin{array}{ll}
F_1(D^2u,\nabla u,u,x)=0 & \text{in }B_1 \cap \Omega ,\\
F_2 (D^2 u,\nabla u,u,x)=0 & \text{in }B_1\setminus\Omega ,\\
u \in \mathbb{S} (B_1),
\end{array}
\right.
\end{equation}
where  $\Omega$ is an unknown open set,   $F_i$ ($i=1,2$) are elliptic/parabolic operators,
and $\mathbb{S}(B_1)$ is some function space to which $u$ should belong.\\

This general setting includes as special cases several important problems, for instance:
\begin{enumerate}
\item[-] The classical no-sign obstacle problem, that includes the classical obstacle problem as  particular case: $F_1:=\Delta u-1$, $F_2:=u$, $\mathbb{S}(B_1):=W^{2,n}(B_1)$
\cite{ALS,PSU}.
\item[-] The general no-sign obstacle problem for fully non-linear operators when $\psi \in C^{1,1}$: $F_1:=F(D^2u)-F(0)$, where $F$ is any convex fully-nonlinear operator, 
$F_2 :=\chi_{\{|D^2 u|> C_0\}}$ (since on the contact set $|D^2u|=|D^2\psi|\leq C_0$), and $\mathbb{S}(B_1)=W^{2,n}(B_1)$ \cite{FS}.
\end{enumerate}
The global $W^{2,n}$ assumption on $u$ is motivated by the theory developed in \cite{CCKS} (see also \cite{F, CS, PSU}) and it is crucial
for regularity,
as it implies that the two equations are not unrelated. 
In both cases one can prove optimal regularity for the solution, namely $u \in C^{1,1}$.
However this optimal regularity result is specific to the structure of these two problems and cannot be true 
in general, as the following examples show.

   \begin{example}\label{unstable-ex}
  In \cite{AW} the authors construct an interesting example of a function $u$ which is not $C^{1,1}$ and 
  solves 
  $$\Delta u = -\chi_{\{u>0\}}  \qquad \hbox{in } B_1.$$
   This function is a solution to the above problem
  with $F_1(D^2u)= \Delta u +1$ in $\Omega=\{u>0\}$, and $F_2(D^2u)= \Delta u $ in $B_1\setminus \Omega$.
  Moreover $\Omega$ is a level set of $u$, but the $C^{1,1}$ regularity fails. 
  Notice that by elliptic regularity $u \in C_{\rm loc}^{1,\alpha}(B_1)$ for any $\alpha <1$, and actually the second derivatives of $u$ belong to BMO.     
  
  For systems, Nina Uraltseva gave the following simple example: $ S(z) := z^2\log |z|$, $z=x+ {\bf i}y$. The real and imaginary parts 
  of $S$ satisfy (up to a multiplicative constant) $ \Delta u_i=\frac{- u_i}{|{\bf u}|}$, $(i=1,2)$, but they are not $C^{1,1}$.  
   \end{example}

The examples above suggest that in general we cannot expect $C^{1,1}$-regularity of solutions, unless we  restrict 
the problem into a smaller class of equations; this is illustrated in a few cases below.

Notice that in Example \ref{unstable-ex} although the solution is not $C^{1,1}$ their second derivatives  belong to BMO  as a consequence of classical elliptic regularity (see, for instance, the appendix in \cite{FS2} for a proof of BMO regularity
for elliptic and parabolic fully-nonlinear equations). Nevertheless, we do not expect this fact to be true in general for solutions to \eqref{matching},
and it
would be interesting to understand 
under which assumptions such BMO regularity holds.\\

We conclude this section by noticing that in the book \cite[Chapter VI, Section 5]{Kind-stam} (see also \cite{KNS}) there are discussions and  few results on ``overdetermined'' problems with matching Cauchy data. 
However, there the free boundary is a priori  assumed to be $C^1$ and the authors derive higher order regularity,
while here we want to start with much less regularity for the solutions and without any regularity assumption  on the free boundary.  

A similar type of problem is the one of ``matching Cauchy data'' from one side of a domain\footnote{
Notice that formally, if $u_1=u_2$ in the complement of $\Omega$, then their gradients agree there and we are saying that
$u_1$ and $u_2$ are functions solving some nice equation inside $\Omega$ and satisfy both $u_1=u_2$ and $\nabla u_1=\nabla u_2$
on $\partial\Omega\cap B_1$.} 
\begin{equation}\label{cauchy}
\left\{
\begin{array}{ll}
F_i(D^2u_i,\nabla u_i,u_i,x)=0 & \text{ in }B_1 \cap \Omega ,    \qquad (i=1,2)\\
u_1=u_2,  & \text{ in }B_1 \setminus \Omega ,\\
u_1, u_2 \in W^{2,n}(B_1),\\
\end{array}
\right.
\end{equation}
with $F_i$ ``nice'' elliptic operators, $F_1(M,p,z,x) \neq F_2(M,p,z,x)$ for all $(M,p,z,x)$ (so that in particular $u_1\not\equiv u_2$ in $\Omega$), and one may ask both about the regularity of solutions
and that  of $\partial \Omega$. A simple example of such a problem in the linear case is the following:
$$
F_1(D^2u_1,\nabla u_1,u_1,x)=\Delta u_1  , \qquad F_2(D^2u_2,\nabla u_2,u_2,x)=\Delta u_2 - 1,
$$ 
in $B_1\cap \Omega$. Then $u=u_2 - u_1$ satisfies $\Delta u= 1$ in $\Omega \cap B_1$, and $u =0$ in
$B_1 \setminus \Omega $. As mentioned before, in this case $u \in C^{1,1}$, and starting from there one can obtain regularity for the free boundary (under some thickness assumptions of the complement of $\Omega$)
as in \cite{CKS}.
However, for the nonlinear problem one needs a different approach as one cannot subtract $u_1$ and $u_2$, and to our knowledge
there are no results in this direction.


 \section{Optimal regularity for unconstrained problems}

\subsection{The double obstacle problem}   \label{sect:double}

The constrained variational problem,  known as double obstacle problem, is given by
\begin{equation}\label{double-obst}
\min_{\mathcal K} \int_{B_1} |\nabla v|^2 \qquad \qquad {\mathcal K}:=\{v \in W_g^{1,2} :\  \psi_1 \leq v\leq \psi_2\}.
\end{equation}
Notice that in the particular case of $\psi_2 \equiv +\infty$
the problem reduces to the classical obstacle problem which is more or less very well understood  \cite{Caff,Caobstacle}. Notwithstanding this,  a huge number of interesting applications  in diverse areas in analysis have arisen in the recent years; among them are  the theory of quadrature domains, and related problems such as Laplacian growth, algebraic droplets of Coulomb gas ensembles, potential theoretic equilibrium measure and its random matrix models (see \cite{TBAZW}).

The double obstacle case, however, is much less studied; this depends on less developed appropriate technical tools, and (probably) on the fact that it has been considered as an ``easy'' generalization of the one-sided obstacle problem. This is of course the case if one considers the problem at points  away from the part of the free boundary that is ''sandwiched''  between both obstacles.

Here, we want to show  how it is possible to restate this problem in a form that allow us to deduce regularity of solutions as a corollary 
of general results for unconstrained problems. As we shall see, it is actually possible to replace the Laplacian by a fully nonlinear operator and still get optimal regularity.

Before passing to fully nonlinear operators, let us make the following observation:
assume that $u$ solves some equation of the form 
$$
F(D^2u,\nabla u,u,x)=0
$$
and suppose that $u \in C^{1,\alpha} (B_1)$ for some $\alpha>0$. Then, if $F$ is at least H\"older continuous in the variables $(\nabla u,u,x)$, the operator
$$
G(M,x):=F\bigl(M,\nabla u(x),u(x),x\bigr)  \qquad \hbox{in } B_1
$$
is independent of $u$ and $\nabla u$, and is H\"older continuous in the $x$ variable. This is to say that when it comes to prove high enough regularity (say $C^{1,1}$)
one can assume that the operator is independent of  $u$ and $\nabla u$, since the above argument allows one easily to reduce to that case.

Thus, we fix
two functions $\psi_1 \leq \psi_2$ of class $C^{1,1}$, and consider the solution of the problem
$$
\max\bigl\{\min \{ -F(D^2u,x) , u-\psi_1\}, u-\psi_2 \bigr\}=0\qquad \text{in $B_1$}
$$
with $\psi_1<u<\psi_2$ on $\partial B_1$. Here $F$ is a fully nonlinear operator
(that one could make  depend also on $u$ and $\nabla u$, if desired)
and this equation is to be understood in terms of viscosity solutions \cite{CCKS}.
By standard theory \cite{CCKS} one can show that any solution belongs to $W^{2,n}$,
and that on the contact set $\{u=\psi_1\}$ (resp. $\{u=\psi_2\}$) the value of $F(D^2u,x)$
is given by replacing $D^2u$ with $D^2\psi_1$ (resp. $D^2\psi_2$). 
Hence, one can rewrite the above equation as
\begin{equation}\label{reform-double}
F(D^2u,x)= F(D^2\psi_1,x)\chi_{\{u=\psi_1\}} + F(D^2\psi_2,x)\chi_{\{u=\psi_2 \neq \psi_1\}} \qquad \hbox{in } B_1,
\end{equation}
along with the constraint $\psi_1 \leq u \leq \psi_2$.
Since the right hand side is bounded, by classical theory for fully nonlinear equations
one deduces that $u \in W^{2,p} (B_1) $ for any $p\in (1,\infty)$ (actually, $D^2u$ belongs to BMO,
see for instance the appendix in \cite{FS2}), hence $u$  is twice differentiable a.e. 
We then notice that, if we set $\Omega:=\{\psi_1<u<\psi_2\}$, $u$ solves
$$
F(D^2u,x)=0 \qquad\text{ in $\Omega \cap B_1$},
$$
 and 
 $$
 |D^2 u|\leq C_0 \qquad \text{in $B_1\setminus \Omega$},
 $$
 (since there $D^2u$
is equal to $D^2\psi_1$ or $D^2\psi_2$).
Hence, as a corollary of the results in \cite{FS} (see also \cite{IM}) we conclude that  $u$ is of class $C^{1,1}$ inside $B_{1/2}$.
More precisely, the following theorem holds:
\begin{theorem}\label{double}
Let $F(M,x)$ be a fully-nonlinear uniformly elliptic operator satisfying
\begin{equation}
\label{eq:F holder}
|F(M,x) -F(M,y)| \lesssim \bigl(1+|M|\bigr)|x-y|^\alpha
\end{equation}
for some $\alpha>0$.
Assume that $F(\cdot,x)$ is either convex or concave, and that $\psi_1,\psi_2 \in C^{1,1}$.
 Then for any solution $u$ to equation \eqref{reform-double} one has $u \in C^{1,1}(B_{1/2})$, where the $C^{1,1}$ norm of $u
$ in $B_{1/2}$ depends only on the data and $\|u\|_{L^\infty(B_1)}$.
\end{theorem}



\subsection{Matching regularity}
A general question (that includes the above discussion)   is the following:\\ 

{\it Let $u:B_1\to {\mathbb{R}}$ belong to some function space $\mathbb{S}(B_1)$ (to be specified). Suppose that $u$ solves 
a PDE in $\Omega \cap B_1$, has some smoothness in $B_1\setminus \Omega$,
and the PDE is of such character that it implies the same smoothness in $\Omega$. Can we expect that $u$ enjoys such a regularity also across $\partial \Omega$?\\
}

Let us formulate this problem in some concrete examples.\\

\begin{example}\label{ex-gradient} Consider the problem
\begin{equation}\label{lower-regul-1}
\left\{
\begin{array}{ll}
F(D^2u,x)=0 & \text{in }B_1 \cap \Omega ,\\
\nabla u  \in C^\alpha & \text{in }B_1\setminus\Omega ,\\
u \in W^{2,n} (B_1),
\end{array}
\right.
\end{equation}
where $\alpha\in (0,1)$ and $F$ is a "nice'' uniformly elliptic fully nonlinear operator.
Then we ask whether $u\in C^{1,\alpha} (B_{1/2})$.
\end{example}
The above question is nontrivial already when $F=\Delta u$. 
Notice that if instead of saying that $\nabla u  \in C^\alpha$ one knows that $D^2u \in L^p(B_1\setminus \Omega)$ for some $p>n$,
then $F(D^2u,x) \in L^p$ in $B_1$ and the $W^{2,p}$ regularity of $u$ in $B_{1/2}$ follows immediately from
elliptic regularity \cite{CaFNL}.

Similar questions can also be asked for more degenerate operators of the form $|\nabla u|^\gamma F(D^2u,x)$
that are known to satisfy interior $C^{1,\beta}$ estimates \cite{IS}.\\

Example \ref{ex-gradient} is just a simple illustration of some questions that naturally
include several free boundary problems as a special case.
For instance, consider again a $W^{2,n}$ solution of the no-sign obstacle problem
\begin{equation}
\label{eq:obst}
 \Delta u \,\chi_{\{u\neq \psi\}}=0. \footnote{So far, this seems to be the most general formulation of obstacle-type problems that encompasses many known free boundary problems.  Notice that this formulation (in the viscosity sense) does not require any a priori regularity of $\psi$ besides continuity.}
\end{equation}
If $\psi \in C^{1,\alpha}$ then the regularity of $u$ can be seen as a particular case of \eqref{lower-regul-1} above
with $\Omega=\{u\neq \psi\}$.
More  generally, one may ask for regularity of solutions under other regularity assumptions on $\psi$.

We think that a better  understanding  these problems   would be extremely interesting not only
for the applications to free boundaries, but also because their study could lead to the development of new interesting
 techniques.\\

A different direction worth pursuing would be to strengthen some of the conditions on $u$ and ask for higher regularity.
We illustrate this with another example.

\begin{example}\label{ex-higher-regul} Consider the solution to
\begin{equation}\label{lower-regul}
\left\{
\begin{array}{ll}
F(D^2u,x)=0 & \text{in }B_1 \cap \Omega ,\\
|F(0,x)| \lesssim |x|^\alpha  & \text{in }B_1 \cap \Omega\\ 
u \in W^{2,n} (B_1).
\end{array}
\right.
\end{equation}
where $\alpha \in (0,1]$, and assume that
\begin{equation}
\label{eq:small D2u}
|D^2 u (x)| \lesssim |x|^\alpha   \qquad \text{a.e. in }B_1\setminus\Omega.
\end{equation}
Then we ask whether $u$ is $C^{2,\beta}$ at the origin for some $\beta>0$, that is, whether there is a second degree polynomial  $P(x)$ such that 
$$
|u(x)-P(x)|\lesssim |x|^{2+\beta} \qquad \forall \,x \in B_1.
$$
\end{example}
In this particular case, the answer to this question is a consequence of the results in \cite{CaFNL}.
More precisely, assuming for simplicity that $F(D^2u,x)=G(D^2u)-f(x)$ with $G(0)=0$, the above assumptions imply that
$$
G(D^2u)=f
$$
with 
$$
|f(x)| \lesssim |x|^\alpha.
$$
Then, if $G$ is uniformly elliptic and either convex or concave, 
\cite[Theorem 3]{CaFNL} implies 
the $C^{2,\beta}$ regularity at the origin for any $\beta<\min\{\alpha,\alpha_0\}$,
where $C^{2,\alpha_0}$ is the regularity for the clean equation $G(D^2u)=0$.

To show a simple case where this result can be applied consider the no-sign obstacle problem for the Laplace operator:
assume that $\psi \in C^{2,\alpha}(B_1)$ and let $u \in W^{2,n}$ solve 
$$
\Delta u\,\chi_{\{u\neq \psi\}}=0.
$$
Then $v:=u-\psi$ is a solution of 
$$
\left\{
\begin{array}{ll}
\Delta v=-\Delta \psi & \text{in }B_1 \cap \Omega ,\\
v=0   & \text{in }B_1 \setminus \Omega ,\\ 
\end{array}
\right.
$$
with $\Omega:=\{v\neq 0\}$. 
Then $v$ is universally  $C^{1,1}$ in $B_{1/2}$ and we can distinguish two kind of points, depending 
whether $\Delta\psi(x_0)\neq 0$ or not.
In the first case we are at the so-called ``non-degenerate points'' which are the ones where regularity of the free boundary can be proved (see for instance \cite{Caobstacle,PSU}). When $\Delta\psi(x_0)=0$ then $|\Delta \psi(x) - \Delta\psi(x_0)| \lesssim |x-x_0|^\alpha$ 
and the above discussion implies that 
$u\in C^{2,\gamma}$ at $x_0$ for all $\gamma<\alpha$.

%
%

\section{Systems and Switching problems}
\subsection{Optimal switching problems}
Optimal switching problems, where a state is switched to another state for cost reduction or profits, corresponds to obstacle/constrained problems where vectorial functions are involved.

Switching problems have recently attracted a lot of attention in mathematical finance and economics, where    uncertainty related to  evaluation of investment projects affects decisions.
  Such problems also  arise in the optimal control of hybrid systems and in
 stochastic switching zero-sum game problems; see for instance \cite{EF79,ADPS,DHP} and references therein. For
applications to starting-stopping problem and finance, we refer the reader to
\cite{HJ07,DH09}.

Assuming as before that our operators do not depend on $u$ and $\nabla u$,
the problem of optimal switching for two-states are formulated in terms of viscosity solutions to the following system of equations (in $B_1$, say):
\begin{equation}\label{opt-switch-1}
\left\{
\begin{array}{ll}
&\min\bigl\{-F_1(D^2u_1,x), u_1 -(u_2 -\psi_1)\bigr\}=0 , \cr
&\min\bigl\{-F_2(D^2u_2,x), u_2 -(u_1 -\psi_2)\bigr\}=0.
\end{array}
\right.
\end{equation}
Here as usual $F_i$ are nice elliptic operators, $\psi_1 + \psi_2 \geq 0$, and both $\psi_1$ and $\psi_2$ are smooth. We may rephrase this, upon showing a lower-regularity of type $W^{2,n}$, in the form
\begin{equation}\label{opt-switch-2}
\left\{
\begin{array}{ll}
&F_1(D^2u_1,x)= F_1(D^2(u_2-\psi_1),x) \chi_{\{u_1= u_2 - \psi_1\}} , \cr
&F_2(D^2u_2,x)= F_2(D^2(u_1-\psi_2),x) \chi_{\{u_2= u_1 - \psi_2\}} ,
\end{array}
\right.
\end{equation}
along with the constraints $u_i \geq u_j -\psi_i$ ($i\neq j$, and $i,j=1,2$).\\

It is actually possible to go further and consider the case of a double switch.
Then, a double-constrained version for two-switching problem can be formulated as
\begin{equation}\label{opt-switch-3}
\left\{
\begin{array}{ll}
&\max\left\{ \min\bigl\{-F_1(D^2u_1,x), u_1 -(u_2 -\psi_1)\bigr), u_1 - (u_2 -\tilde \psi_1) \right\}=0, \cr
&\max\left\{\min\bigl\{-F_2(D^2u_2,x), u_2 -(u_1 -\psi_2)\bigr),  u_2 - (u_1 -\tilde \psi_2) \right\}=0,
\end{array}
\right.
\end{equation}
with conditions 
\begin{equation}\label{constrained-2}
\tilde \psi_i \leq \psi_i, \qquad \psi_1 + \psi_2 \geq 0, \qquad \tilde \psi_1 + \tilde \psi_2 \geq 0.
\end{equation}
We refer to \cite{DHMZ} and the references therein for some background on this family of problems.

As for \eqref{opt-switch-1} it is not hard to realize that, once the ``primary''  (here it is $W^{2,n}$) regularity for solutions is established
 one may rewrite \eqref{opt-switch-3} in the form
\begin{equation}\label{opt-switch-3bis}
\left\{
\begin{array}{ll}
F_1(D^2u_1,x)=   &F_1(D^2(u_2-\psi_1),x) \chi_{\{u_1= u_2 - \psi_1\}} \\
&+  F_1(D^2(u_2-\tilde \psi_1),x) \chi_{\{u_1= u_2 - \tilde \psi_1, \psi_1\neq \tilde \psi_1\}},\cr
 F_2(D^2u_2,x)= &F_2(D^2(u_1-\psi_2),x) \chi_{\{u_2= u_1 - \psi_2\}}\\
&  +  F_2(D^2(u_1-\tilde \psi_2),x) \chi_{\{u_1= u_2 - \tilde \psi_1,  \psi_2\neq \tilde \psi_2\}} ,
\end{array}
\right.
\end{equation}
with the constraints \eqref{constrained-2} for $u_1,u_2$.

\subsection{Regularity issues for unconstrained system}
To define an unconstrained version of  \eqref{opt-switch-3bis}, it suffices to drop condition  \eqref{constrained-2}. Nevertheless, we shall rephrase equations \eqref{opt-switch-3bis} to allow further generalization of an unconstrained problem. We  thus define 
$$
A_i:=\{u_i = u_j - \psi_i\},\quad \tilde A_i=\{u_i = u_j - \tilde \psi_i\}, \qquad \text{for $i,j=1,2$ and $i\neq j$},
$$
and set 
$$
\Omega_1 := B_1\setminus (A_1 \cup \tilde A_1),\qquad \Omega_2 := B_1\setminus (A_2 \cup \tilde A_2).
$$
This allows us to rephrase \eqref{opt-switch-3bis} in the following general form:
\begin{equation}\label{unconstrained-switch-2}
\left\{
\begin{array}{ll}
F_1(D^2u_1,x)= 0        &\text{in } B_1\cap  \Omega_1     ,\cr
F_2(D^2u_2,x)= 0        & \text{in }B_1\cap \Omega_2     ,\cr
  |D^2 (u_1 - u_2)| \leq C & \text{in } B_1\setminus (\Omega_1\cup\Omega_2).
\end{array}
\right.
\end{equation}
Then, assuming that
$\hbox{dist}(B_1\setminus \Omega_1, B_1\setminus \Omega_2)=\delta >0$
(this is satisfied for instance if we impose that $\psi_1+\tilde \psi_2>0$ and $\psi_2+\tilde \psi_1>0$)
we can immediately prove  uniform $C^{1,1}$-regularity for $u_1$ and $u_2$ in $B_{1/2}$. 

Indeed, if $N_2$ denotes the  $\delta/4$-neighborhood of $B_1\setminus \Omega_2$, then $u_2$ is uniformly  $C^{2,\alpha}$ in  $B_{3/4}\setminus N_2$ because it solves a nice PDE in a $\delta/2$-neighborhood of this set, with bounded boundary data. In particular, in $B_{3/4}\setminus N_2$, $u_1$ is a solution to a free boundary problem of the type studied in \cite{FS,IM}, hence $u_1$ is $C^{1,1}_{\rm loc}$ there. Since $u_1$ also solves a nice PDE 
in a $\delta/2$-neighborhood of $N_2$, we conclude that $u_1$ is uniformly $C^{1,1}$ in $B_{1/2}$. 
The same argument applies to $u_2$. For future reference, we formulate this as a theorem.


\begin{theorem}\label{switch}
Let $F_i(M,x)$ satisfy condition \eqref{eq:F holder}, and assume that
$F_i(\cdot,x)$ is either convex or concave, and that $\psi_1,\psi_2,\tilde \psi_1,\tilde\psi_2 \in C^{1,1}(B_1)$.
Also, let $u_1, u_2$ solve \eqref{unconstrained-switch-2} and assume that $\hbox{dist}(B_1\setminus \Omega_1, B_1\setminus \Omega_2)=:\delta >0$. Then $u_1, u_2  \in C^{1,1}(B_{1/2})$.
\end{theorem}

Let us point out that the situation $\delta=0$ is much more complicated, and already the case with one switch is far from trivial.
Very recently, in \cite{A} the author has obtained the optimal $C^{1,1}$-regularity for the minimal solution
of \eqref{opt-switch-1} when $F_1=\Delta u -f_1$ and $F_2=\Delta u-f_2,$ under the assumption that the zero loop set $\{\psi_1+\psi_2=0\}$ is the closure of its interior. As shown in the same paper, this result is optimal,
and it is possible to find a counterexample showing that the $C^{1,1}$-regularity does not hold without that assumption on the zero loop.

It would be extremely interesting to understand whether the same regularity result can be shown for the unconstrained case, that is, when one forgets the constraints
$u_1-u_2 +\psi_1 \geq 0$ and $u_2-u_1 +\psi_2 \geq 0$.

%
%

\section{Diversifications and open questions}
In this section we shall present  several diverse but related free boundary problems, where many of them  so far have been almost untouched. We shall formulate the problems without entering into much details, and suggest new problems that may be of interest for future research. 

\subsection{Gradient constraints}

A classical problem is the well-known ``gradient constrained problem'', that was subject to intense study 
a few decades ago. In its simplest form, it reads
\begin{equation}\label{gradient-obst-1}
\min_{\mathcal K} \int_{B_1} |\nabla v|^2 \qquad \qquad {\mathcal K}:=\{v \in W_g^{1,2} :
\  |\nabla v|\leq h(x) \}.
\end{equation}
This problem has applications in elastoplasticity of materials \cite{Tsu} and optimal control problems \cite{Son-Shr}. There are also recent applications
in mathematical finance, where transactions costs are involved, see  \cite{Son-Shr2} and the references therein.

When $h$ is smooth one can show that solutions are $C^{1,1}$, see  \cite{W}.
However, without the smoothness of $h$ the problem becomes much more delicate. 
A possible way to attack this problem in the case $h=1$ (or more in general when $\Delta(h^2)\leq 0$, as indicated in \cite[Theorem 3.1]{Santos}) is to rewrite the problem as a double obstacle problem, where the obstacles solve the Hamilton-Jacobi equation $\pm |\nabla u|=\pm h$ in the viscosity sense. 

This motivates the study of double obstacle problems where the obstacles solve a general Hamilton-Jacobi equation.
Since in general solutions to the Hamilton-Jacobi equation
are not smooth, the regularity of solutions to this double obstacle problem is far from trivial.
We refer to the recent article  \cite{ASW}, and the references therein.

In \cite{ASW} it was shown that if the constraint is set as $|\nabla u -a(x)| = 1$ with $a$ a vector field of class $C^\alpha$,
and one assumes that $u$ is $C^{1}$ near the contact region, then one obtains $C^{1,\alpha/2}$ regularity for the Hamilton-Jacobi equation
 and from there one may proceed  to obtain  $C^{1,\alpha/2}$ regularity for the solutions.

Motivated by the discussion above,
we define now a  general  class of (double) gradient constrained problems for fully nonlinear equations
as 
$$
\max\bigl\{\min \{ -F(D^2u,x) , |\nabla u -a|-h_1\}, |\nabla u-a|-h_2 \bigr\}=0\qquad \text{in $B_1$},
$$
where the equation is in the viscosity sense.  Formally, the equation can  be written as 
$$
F(D^2u,x)\chi_{\{|\nabla u-a| \neq h_1 \} \cup \{|\nabla  u-a |\neq h_2\} }=0,\qquad
h_1 \leq |\nabla u-a| \leq h_2.
$$
Dropping the constraints, one encounters a completely new  unconstrained problem.\\
A natural question is: {\it How regular are  solutions and the free boundary for such problems?}

\subsection{The $p$-Laplace operator}
The case of the $p$-Laplace operator, i.e.  $\Delta_p u:= \hbox{div} (|\nabla u|^{p-2}\nabla u)$ ($1<p<\infty$), introduces challenging difficulties, and one may find a large number of results in the literature concerning  lower-order regularity of solutions.
While solutions to the $p$-Laplacian are no better than $C^{1,\alpha_p}$ for some (unknown) exponent $\alpha_p \in (0,1)$,
it is interesting to observe the regularity improves near the free boundary: indeed, it has been shown recently in \cite{ALS2} that solutions are $C^{1,1}$ {\it at free boundary points}.\footnote{This means that the maximal difference between the solution and the obstacle,  $u-\psi$, on a ball $B_r(z)$, with $z$ on the free boundary, 
is controlled by $r^2$.} 

The unconstrained counterpart of this problem, that is
$$
\Delta_p u =\Delta_p \psi \chi_{\{u=\psi\}} 
$$
without any restriction of the type $u\geq \psi$, is a completely untouched area. Low regularity can be obtained up to a certain order, due to boundedness of the $p$-Laplacian of $u$. In general this gives (through a blow-up and Liouville theorem) regularity of order 
$1+\alpha$ for any $\alpha < \alpha_p$.
It seems plausible to expect a second order regularity at free boundary points also for the above unconstrained problem, i.e., without the condition $u\geq \psi$. 

Similar formulations can be made for the gradient constrained or unconstrained problems for the case of $p$-Laplacian. One expects that,
also in these cases, the optimal regularity of solutions at  free boundary points should be of order two.

\subsection{Monotone operators in the $u$ variable}
The semilinear problem given by  $F(D^2 u, u):= \Delta u - f(u)=0$
where $f(u)$ is monotone-increasing and has a jump discontinuity for some values of $u$, can be seen as an unconstrained free boundary problem, where free boundary is given by the level surfaces for $u$ where $f(u)$ has a discontinuity.  It was shown in \cite{Sh} that solutions to this problem are $C^{1,1}$. It is tantalizing to analyze the case of fully nonlinear equations
$F(D^2 u, u)=0$ with $F'_u \leq 0 $.
Notice that in the particular case $F(D^2 u, u)=G(D^2u)-f(u)$ with $G$ convex and $f$ bounded, 
elliptic regularity gives that $D^2u$ belongs to BMO, so one only needs to understand the ``last step'' from BMO to $L^\infty$.

 \subsection{Optimal regularity at free boundary points}
Degenerate/singular  PDEs in general fail to have good regularity estimates. E.g., it is well-known that the $p$-harmonic functions in general cannot be better than $C^{1,\alpha_p}$ for a universal exponent $0< \alpha_p <1$ (the exact value of $\alpha_p$ is currently unknown).
 Notwithstanding this, the authors in \cite{ALS2} proved that solutions  to the obstacle problem 
 enjoy second order regularity at the free boundary. Similar results seem plausible for the gradient constrained problems, as well as 
for degenerate/singular equations (e.g. $|\nabla u|^\gamma F(D^2u)$ or $u^\alpha F(D^2u)$).

An even more interesting question is whether this higher regularity can propagate from the free boundary into a neighborhood of it. More precisely:\\
{\it   Is it  true that solutions to such degenerate problems are $C^{1,1}$ in a uniform neighborhood of the free boundary?}

 \subsection{Parabolic problems}
All the problems mentioned in this paper have a natural parabolic counterpart.
In their simplest form, one can replace the operators
$F(D^2u,\nabla u,u,x)$ with $F(D^2u,\nabla u,u,x,t)-\partial_t u$,
and ask analogous questions in this setting.
Notice that optimal regularity  (i.e., $C_x^{1,1}\cap C_t^{0,1} $) for the unconstrained obstacle problem with fully nonlinear operators 
has been recently shown in \cite{FS2,IM}, and we expect that 
results valid in the elliptic setting should carry on also to the parabolic
case.





\begin{thebibliography}{9}



\bibitem{A}
Aleksanyan, G.;
Optimal regularity in the optimal switching problem.
Preprint, 2014.

\bibitem{AW}
Andersson, J.; Weiss, G. S.;
Cross-shaped and degenerate singularities in an unstable elliptic free boundary problem.
{\em J. Differential Equations} 228 (2006), no. 2, 633-640. 

\bibitem{ASW} 
 Andersson, J., Shahgholian, H., Weiss, G. S., Double obstacle problems with obstacles given by non-$C^2$ Hamilton-Jacobi equations. Arch. Ration. Mech. Anal. 206 (2012), no. 3, 779-819. 




\bibitem{ALS}
Andersson J.;  Lindgren E.;  Shahgholian H.;
Optimal regularity for the no-sign obstacle problem. {\it Comm. Pure Appl. Math.}.

\bibitem{ALS2}
Andersson J.;  Lindgren E.;  Shahgholian H.;
Optimal Regularity for Degenerate Obstacle Problems (preprint).

\bibitem{ADPS}
Arnarson, T.; Djehiche, B.; Poghosyan, M.; Shahgholian, H.;
A PDE approach to regularity of solutions to finite horizon optimal switching problems.
{\it Nonlinear Anal. } 71 (2009), no. 12, 6054-6067. 




\bibitem{Caff} Caffarelli L. A.; The regularity of free boundaries in higher dimensions. {\it Acta Math.} 139 (1977), no. 3-4, 155-184. 


\bibitem{CaFNL}
Caffarelli, L. A.; Interior a priori estimates for solutions of fully nonlinear equations.
{\it Ann. of Math. (2)} 130 (1989), no. 1, 189-213.

\bibitem{Caobstacle}
Caffarelli, L. A.;
The obstacle problem revisited.
{\it J. Fourier Anal. Appl.} 4 (1998), no. 4-5, 383-402. 



\bibitem{CCKS}
Caffarelli, L.; Crandall, M. G.; Kocan, M.; Swiech, A.; On viscosity solutions of fully nonlinear equations with measurable ingredients. {\it Comm. Pure Appl. Math.} 49 (1996), no. 4, 365-397.




\bibitem{CKS} Caffarelli L. A.; Karp L.; Shahgholian H.; Regularity of a free boundary with application to the Pompeiu problem. {\it Ann. of Math. (2)} 151 (2000), no. 1, 269-292.

\bibitem{CS} Caffarelli, L.; Salazar, J.; Solutions of fully nonlinear elliptic equations with patches of zero gradient: existence, regularity and convexity of level curves. {\it Trans. Amer. Math. Soc.} 354 (2002), no. 8, 3095-3115. 

\bibitem{CSbook}
Caffarelli, L.; Salsa, S.; A geometric approach to free boundary problems.  Graduate Studies in Mathematics, 68. American Mathematical Society, Providence, RI, 2005. x+270 pp.

\bibitem{Choe}
Choe, H. J.; Shim, Y.-S.; Degenerate variational inequalities with gradient constraints. {\it Ann. Scuola Norm. Sup. Pisa Cl. Sci. (4)} 22 (1995), no. 1, 25-53.

\bibitem{DH09}
Djehiche, B.; Hamad\`ene, S.;
On a finite horizon starting and stopping problem with risk of abandonment.
{\it Int. J. Theor. Appl. Finance} 12 (2009), no. 4, 523-543. 

\bibitem{DHMZ}
Djehiche, B.; Hamad\`ene, S.; Morlais, M.-A.; Zhao, X.;
On the Equality of Solutions of Max-Min and Min-Max Systems of Variational Inequalities with Interconnected Bilateral Obstacles.
Preprint, 2014.


\bibitem{DHP}
Djehiche, B.; Hamad\`ene, S.; Popier, A.;
A finite horizon optimal multiple switching problem. 
{\it SIAM J. Control Optim.} 48 (2009), no. 4, 2751-2770. 

\bibitem{TBAZW}
Teodorescu, R.,  Bettelheim, E., Agam, O., Zabrodin, A.; Wiegmann, P.,
Normal random matrix ensemble as a growth problem.  Nuclear Phys. B 704 (2005), no. 3, 407-444. 


\bibitem{EF79}
Evans, L. C.; Friedman, A.;
Optimal stochastic switching and the Dirichlet problem for the Bellman equation. {\it Trans. Amer. Math. Soc.}
253 (1979), 365-389.  


\bibitem{FS}
Figalli, A.; Shahgholian, H.; A general class of free boundary problems for fully nonlinear elliptic equations.
{\it Arch. Ration. Mech. Anal.} 213 (2014), no. 1, 269-286.

\bibitem{FS2}
Figalli, A.; Shahgholian, H.; 
A general class of free boundary problems for fully nonlinear parabolic equations.
{\it  Ann. Mat. Pura Appl.}, to appear.

\bibitem{F} Friedman, A.; Variational principles and free-boundary problems. A Wiley-Interscience Publication. Pure and Applied Mathematics. John Wiley \& Sons, Inc., New York, 1982.

\bibitem{HJ07}
Hamad\`ene, S.; Jeanblanc, M.;
On the starting and stopping problem: application in reversible investments. 
{\it Math. Oper. Res.} 32 (2007), no. 1, 182-192. 

\bibitem{IS}
Imbert, C.; Silvestre, L.; $C^{1,\alpha}$
regularity of solutions of some degenerate fully non-linear elliptic equations. {\it Adv. Math.} 233 (2013), 196-206.

\bibitem{IM}  Indrei, E.;    Minne, A.; Regularity of solutions to fully nonlinear elliptic and parabolic free boundary problems.
Preprint.

\bibitem{KNS}
Kinderlehrer, D.; Nirenberg, L.; Spruck, J.;
Regularit\'e dans les probl\'emes elliptiques |'a frontiere libre. (French. English summary) 
{\it C. R. Acad. Sci. Paris S\'er. A-B} 286 (1978), no. 24, A1187-A1190. 


\bibitem{Kind-stam}
Kinderlehrer, D.; Stampacchia, G.; {\em An introduction to variational inequalities and their applications.} Pure and Applied Mathematics, 88. Academic Press, Inc., New York-London, 1980.







\bibitem{PSU} Petrosyan, A.; Shahgholian, H.; Uraltseva, N., Regularity of free boundaries in obstacle-type problems. Graduate Studies in Mathematics, 136. American Mathematical Society, Providence, RI, 2012. x+221 pp. ISBN: 978-0-8218-8794-3.




\bibitem{Santos} Santos, L.; Variational problems with non-constant gradient constraints. {\it Port. Math. (N.S.)} 59 (2002), no. 2, 205-248.

\bibitem{Sh} Shahgholian, H.; $C^{1,1}$ regularity in semilinear elliptic problems. {\it Comm. Pure Appl. Math.} 56 (2003), no. 2, 278-281. 






\bibitem{Son-Shr} Shreve, S. E., Soner, H. M.;
A free boundary problem related to singular stochastic control. Applied stochastic analysis (London, 1989), 265-301,  Stochastics Monogr., 5, Gordon and Breach, New York, 1991. 

\bibitem{Son-Shr2}
Shreve, S. E., Soner, H. M.;
Optimal investment and consumption with transaction costs. (English summary) 
Ann. Appl. Probab. 4 (1994), no. 3, 609-692. 

%


\bibitem{Tsu}
Tsuan Wu T.; Elastic-plasitc torsion of a square bar. {\it Trans. Amer. Math. Soc.}
123 (1966), 369-401.


\bibitem{W} 
Wiegner, M.;
The $C^{1,1}$-character of solutions of second order elliptic equations with gradient constraint.
{\it Comm. Partial Differential Equations} 6 (1981), no. 3, 361-371. 
 

\end{thebibliography}
\end{document}